\documentclass[11pt]{amsart}
\usepackage[T1]{fontenc}
\usepackage{lmodern}
\usepackage[utf8]{inputenc}
\usepackage[english]{babel}
\usepackage{amsmath,amssymb,mathrsfs,mathtools,enumitem}
\usepackage[margin=1.1in]{geometry}
\usepackage{microtype}

\usepackage{tikz-cd}
\usepackage{tikz}
\usepackage[hidelinks]{hyperref}
\hypersetup{pdftitle={Calabi--Yau fibrations and  orbifold Mordell}}
\setlist[enumerate,1]{label=(\roman*),ref=\roman*,leftmargin=*,itemsep=3pt}

\newtheorem{theorem}{Theorem}[section]
\newtheorem{lemma}[theorem]{Lemma}
\newtheorem{proposition}[theorem]{Proposition}

\newtheorem{conjecture}[theorem]{Conjecture}
\newtheorem{alphtheorem}{Theorem}

\theoremstyle{definition}

\newtheorem{remark}[theorem]{Remark}

\newtheorem*{con}{Conventions}
\numberwithin{equation}{section}

\newcommand{\QQ}{\mathbb Q}
\newcommand{\CC}{\mathbb C}

\newcommand{\PP}{\mathbb P}

\DeclareMathOperator{\Spec}{Spec}

\DeclareMathOperator{\ord}{ord}

\DeclarePairedDelimiter\powerseries{[\![}{]\!]}
\DeclarePairedDelimiter\laurent{(\!(}{)\!)}

\title[Calabi--Yau fibrations and orbifold Mordell]{Calabi--Yau fibrations and orbifold Mordell}

\author{Ariyan Javanpeykar}
\address{Ariyan Javanpeykar,
IMAPP Radboud University Nijmegen,
PO Box 9010, 6500GL,
Nijmegen, The Netherlands}
\email{ariyan.javanpeykar@ru.nl}

\subjclass[2010]
{14D99 
(14E22, 
14G99, 
11G99)} 

\keywords{Calabi--Yau hypersurface,  Campana's special varieties,  abc conjecture, rational points}

\begin{document}

\begin{abstract}
We construct isotrivial projective families of Calabi--Yau varieties with arbitrarily large  inf-multiplicities and no divisible  fibres over the projective line.    As an application, we show that any single  instance of Campana's Orbifold Mordell Conjecture for a general-type C-pair of genus at most one (with only finite multiplicities) implies that the Weakly Special Conjecture is false.
\end{abstract}

\maketitle
\thispagestyle{empty}

 \section{Introduction}\label{section:introduction}
In \cite{BCJW} it is shown that the Weakly Special Conjecture contradicts Campana's Orbifold Mordell Conjecture for C-pair curves $(C,\Delta)$. More precisely it follows from \cite[Remark~5.6]{BCJW} that the Weakly
Special Conjecture would be false if Campana's Orbifold Mordell
Conjecture held for a single orbifold curve of the form
\[
  \left(\mathbb P^1,\tfrac12[a_1]+\cdots+\tfrac12[a_n]\right),
  \qquad n\geq 5,
\]
with the $a_i$ distinct.  In this short note we refine this result. Namely, we show that \emph{every}
  (finite-multiplicities) instance of Orbifold Mordell in genus zero or one
suffices (no matter how large the multiplicities   are). For example, it would suffice to prove Orbifold Mordell for
$\mathbb P^1$ with multiplicities $(7,7,7)$ at $0,1,\infty$, or for an
elliptic curve with a single marked point of multiplicity $m>1$.

A smooth projective variety $X$ over a field $K$ of characteristic zero
is \emph{weakly special} if no finite \'etale cover of $X_{\overline K}$
admits a dominant rational map to a positive-dimensional smooth projective variety of general type.  Over a number field, Lang's conjecture and the Chevalley--Weil
theorem imply that every smooth projective variety with potentially dense
rational points is weakly special \cite{Langconj,CTSSD}. Harris and
Tschinkel formulated the converse,  following suggestions by Abramovich and
Colliot-Th\'el\`ene \cite[Conjecture~1.2]{HarrisTschinkel}.

\begin{conjecture}[Weakly Special Conjecture]\label{conjecture:wsc}
Let $X$ be a smooth projective weakly special variety over a number
field $K$. Then $X(L)$ is Zariski dense for some finite extension $L/K$.
\end{conjecture}

Campana introduced the class of special varieties and conjectured that
a smooth projective variety over a number field is special if and only
if its rational points are potentially dense
\cite[Conjecture~9.20]{Ca04}. Bogomolov and Tschinkel constructed (the
first) examples of weakly special threefolds that are not special
\cite{BT}. For these threefolds, the Weakly Special Conjecture predicts
potential density, whereas Campana's conjecture predicts that the
rational points are not potentially dense. 

Further examples of weakly
special nonspecial varieties were constructed in \cite{RTW,GRTW} building on the construction of Bogomolov and Tschinkel, and in \cite{BCJW} building on Lafon's construction of an Enriques surface $S$ over $\mathbb{Q}(t)$ with no $\mathbb{C}\laurent{t}$-point \cite{Laf04}. 

The projective examples constructed in this paper are families of Calabi--Yau hypersurfaces over a curve, and are among the simplest to
write down;   the general geometric fibre of the  family is a     diagonal (Fermat) hypersurface.   Our starting observation was that a   nonsingular projective model $X\to  \mathbb{P}^1_{\mathbb{Q}}$ of the degree 6 Calabi--Yau hypersurface    
\[
\{x_0^6 + tx_1^6 + t^2 x_2^6 + t^3 x_3^6 + t^4 x_4^6 + t^5 x_5^6 =0\}  \subset \mathbb{P}^5_{\mathbb{Q}(t)}
\] has an inf-multiple fibre over $t=0$ and no divisible fibres; constructing weakly  special nonspecial varieties from such a family can then be achieved by simply pulling back along suitable endomorphisms of $\mathbb{P}^1$. 

Following
\cite{Ca05,Ca11}, let $(C,\Delta)$ be a \emph{C-pair} with  $C$ a smooth projective curve
over a number field $K$ and
$\Delta=\sum_i(1-1/m_i)D_i$ a $\mathbb{Q}$-divisor, where the $D_i$ are distinct closed points
and $m_i\in\mathbb Z_{\geq2}$.  We say that $(C,\Delta)$ is of \emph{general type} if $\deg(K_C+ \Delta)>0$.

We now state Campana's Orbifold Mordell Conjecture (see \cite[Conjecture~4.5 and Remark~4.6]{Ca05}). Let $S$ be a finite set of finite places of $K$ and let $\mathcal{C}$ be a proper regular model over the ring $\mathcal{O}_{K,S}$ of $S$-integers of $K$. Let $\mathcal{D}_i$ denote the closure of $D_i$ in $\mathcal{C}$. For a finite extension $L/K$ and a
finite set $T$ of finite places of $L$ containing all the places of $L$ lying over $S$, a \emph{Campana $\mathcal{O}_{L,T}$-point of $(C,\Delta)$} is a point $P\in C(L)$ such that the associated $\mathcal{O}_{L,T}$-point $\sigma\colon \Spec \mathcal{O}_{L,T}\to \mathcal{C}$ of $\mathcal{C}$
satisfies the following condition. For each $i$ with $P\notin D_i$
and each finite place $v\notin T$, the intersection multiplicity of
$\sigma$ with $\mathcal{D}_i$ at $v$ is either zero or at least $m_i$.
\begin{conjecture}[Orbifold Mordell]\label{conj:orb_mor}
If $\deg(K_C+\Delta)>0$, then $(C,\Delta)$ has only finitely many
Campana $\mathcal{O}_{L,T}$-points for every such $L$ and $T$.
\end{conjecture}

When $\deg(K_C+\Delta)\leq0$, Campana points are potentially dense
\cite{Ca05}. If $\deg(K_C+\Delta)>0$ and $C$ has genus at least two, then Orbifold Mordell follows from
Faltings's theorem \cite{Faltings}. It remains open for every general type C-pair of
genus at most one (with only finite multiplicities).

\begin{alphtheorem}\label{theorem:A}
Suppose that Orbifold Mordell holds for some smooth proper C-pair
$(C,\Delta)$ of general type over a number field with $g(C)\leq1$. Then the Weakly Special Conjecture is
false.
\end{alphtheorem}

The $abc$ conjecture over number fields implies Orbifold Mordell
\cite[Theorem~5.3]{Smeets}, building on
\cite{Elkies,AbramovichBirGeom,Ca05}. As van Schaick shows in   \cite{vanSchaick},  a much weaker ``polynomial'' version of the abc conjecture suffices to prove Campana's Orbifold Mordell Conjecture in certain cases. Namely, van Schaick shows that such a polynomial version of abc implies Orbifold Mordell for $\mathbb P^1$ with sufficiently
large multiplicities at $0,1,\infty$. Thus
\[
\begin{aligned}
  \text{Polynomial $abc$ over number fields}
    &\xRightarrow{\text{van Schaick}}
      \text{one such instance of Orbifold Mordell},\\[3pt]
  \text{one such instance of Orbifold Mordell}
    &\xRightarrow{\text{Theorem~\ref{theorem:A}}}
      \text{the Weakly Special Conjecture is false}.
\end{aligned}
\]

To prove Theorem~\ref{theorem:A}, we construct a weakly special variety
fibred over $C$ whose orbifold base divisor dominates $\Delta$.

For a fibration $f\colon X\to C$ over $K$ and a closed point
$c\in C$, write $f^*(c)=\sum_j a_{c,j}F_{c,j}$. Its
\emph{inf-multiplicity} is $m_c=\min_j a_{c,j}$, and its
\emph{gcd-multiplicity} is $\gcd_j a_{c,j}$. The former determines the
orbifold base $(C,\Delta_f)$, where $\Delta_f=\sum_c(1-1/m_c)c$.
The fibre is \emph{divisible} if its gcd-multiplicity is greater than one.
For example, a fibre of the form $2F_1+3F_2$ is inf-multiple
(its inf-multiplicity is two), but its gcd-multiplicity is one.

\begin{alphtheorem}\label{theorem:B}
Let $(C,\Delta)$ be a smooth proper C-pair of general type over a
number field $K$, with $\Delta = \sum_{i} \left(1-\frac{1}{m_i}\right) D_i$ and  $g(C)\leq1$.
Let $S$ be the reduced support of $\Delta$, and choose primes $p<q$
with $p\geq m_i$ for every $i$. There exist a smooth projective weakly
special variety $X$ over $K$ and a fibration $f\colon X\to C$ such that
\begin{enumerate}
\item $f$ is isotrivial and its generic fibre is a smooth hypersurface
      of degree $pq$ in $\mathbb P^{pq-1}_{K(C)}$;
\item $f$  has orbifold base $(C,(1-1/p)S)$ and has no
      divisible geometric fibres.
\end{enumerate}
If Orbifold Mordell holds for $(C,\Delta)$, then $f(X(L))$ is finite
for every finite extension $L/K$.
\end{alphtheorem}
   With notation as in Theorem \ref{theorem:B},  Campana's key lemma on integral points implies that  $L$-points of $X$ project to Campana $\mathcal{O}_{L,T}$-points of
$(C,\Delta_f)$ for some $T$,  and hence of $(C,\Delta)$
by the inequality $\Delta_f=(1-1/p)S\geq\Delta$ (we refer the reader to the proof of Theorem \ref{theorem:B} for details).  Thus $X(L)$ lies in finitely
many fibres (and is thus not dense in $X$),  assuming $(C,\Delta)$ satisfies Orbifold Mordell. Also, since  $\Delta_f \geq \Delta$, the orbifold base of $f$ is of general type, so that   $X$ is not special. 

The base curve $C$ has genus at most one and is thus weakly special.
Every smooth geometric fibre of $f$ has trivial canonical bundle and is
therefore weakly special. Since $f$ has no divisible geometric fibres
by~(ii),  the fibration criterion \cite[Theorem~2.3]{BCJW} implies that
$X$ is weakly special.

The construction of $X$ in Theorem \ref{theorem:B} starts with Lang's diagonal hypersurface
\cite[p.~377]{Lang52}
\begin{equation}\label{equation:intro_diagonal}
  V_{pq}\colon
  x_0^{pq}+t x_1^{pq}+\cdots+t^{pq-1}x_{pq-1}^{pq}=0
  \qquad\text{over }\mathbb Q(t).
\end{equation}
 Lang showed that this hypersurface has no $\mathbb{C}\laurent{t}$-points.  In fact,   $V_{pq}$ has a point over
$\mathbb C\laurent{t^{1/e}}$ exactly when $p\mid e$ or $q\mid e$.  Thus, its
index over $\mathbb C\laurent{t}$ is one, while the least degree of a closed
point is $p$ (see also   \cite[Proposition~4.4]{ELW} for the statement about the index).
The base change $t=z^{q-p}(z-1)^p$ leaves a single inf-multiple
geometric fibre, over $z=0$, since the orders at $0,1,\infty$ are
$q-p,p,-q$, respectively. A further pullback places this fibre over
the prescribed support on $C$.

Campana and Stoppino constructed non-divisible inf-multiple fibres in
families of higher-genus curves \cite{Ca05,Stoppino}. For projective fibrations in
abelian varieties over a smooth complex curve, the two multiplicities agree
\cite[Corollary~2.2]{CKV}.  Non-divisible fibres of inf-multiplicity two in projective families
of K3 surfaces over $\PP^1$ were first constructed in \cite[Theorem~C]{BCJW}. Our Theorem~\ref{theorem:B} shows that for
simply connected Calabi--Yau fibres the inf-multiplicity can be an
arbitrarily large prime while the gcd-multiplicity is one. For
restrictions on divisible fibres when the total space is itself simply
connected with trivial canonical bundle, see \cite[Theorem~1.3]{KOS}.

\begin{con}
A variety over a field is a geometrically integral separated scheme of finite type.
\end{con}

\subsection*{Acknowledgements} This paper grew out of the joint works \cite{BCJW,BJCampanaStacks}.  The example in   this paper was found whilst preparing for a course given in June 2026 at the Ecole Polytechnique organized by the ANR project AdAnAr; I am  grateful to Carlo Gasbarri and Gerard Freixas i Montplet for the invitation.    I am grateful to Finn Bartsch for the many inspiring and helpful discussions on Campana's special varieties and C-pairs. I am also grateful to  Fr\'ed\'eric Campana, Erwan Rousseau, and Olivier Wittenberg for  the many fruitful discussions over the past years.
This work was supported by the Netherlands Organization for Scientific Research NWO, through grant OCENW.M.23.206.

\subsection*{Use of AI}
ChatGPT 5.6 was used to proofread the paper, revise the introduction,
and double-check proofs and statements.   The example in this paper was first considered by the author for $d=2\cdot 3$.  The idea to consider $d=p\cdot q$ was   discussed with ChatGPT.

\section{Families of Calabi--Yau hypersurfaces}\label{section:geometry}

The following lemma is stated over $\CC\laurent{t}$ in
\cite[Lemma~3.6]{BCJW}. The same proof applies over $k\laurent{t}$ for
every algebraically closed field $k$ of characteristic zero:
each component multiplicity occurs as the degree of a closed point
by \cite[Lemma~2.4]{KestelootNicaise}, and intersecting a multisection
with the special fibre expresses its degree as a nonnegative integral
linear combination of the component multiplicities.
 To state the lemma, recall that the index $i(V)$ of a variety $V$ over a field $k$ is the gcd of all integers $n$ such that $V$ has a closed point of degree $n$.

\begin{lemma}\label{lemma:local_index}
Let $k$ be an algebraically closed field of characteristic zero, put $K=k\laurent{t}$ and $K_e=k\laurent{t^{1/e}}$, and let $V$ be a smooth proper variety over $K$. Let $\mathcal V$ be an integral regular proper flat model over $k\powerseries{t}$, and write $
 \mathcal V_0=\sum_{j=1}^r a_jF_j.
$
Then
\[
 \min_j a_j=\min\{e\geq1:V(K_e)\neq\emptyset\},
 \quad \textrm{and} \quad 
 \gcd_j a_j=i(V).
\]

\end{lemma}

\subsection{Lang's diagonal Calabi--Yau hypersurfaces}\label{subsection:diagonal}

For an integer $d\geq4$ and a field $k$ of characteristic zero, let
\begin{equation}\label{equation:diagonal}
 V_d=Z\left(\sum_{i=0}^{d-1}t^ix_i^d\right)
 \subset\PP^{d-1}_{k(t)}.
\end{equation}
The defining polynomial is a special case of Lang's construction of normic forms \cite[p.~377]{Lang52}.  Lang's observation that $V_d$ has no $\overline{k}\laurent{t}$-point
implies, by Lemma~\ref{lemma:local_index} after extending the constants
to $\overline{k}$, that the geometric fibre at $t=0$ of every regular
projective model over $k[[t]]$ is inf-multiple. The following proposition gives a more precise control of the closed points of $V_d$.

\begin{proposition}\label{proposition:diagonal}
The following statements hold.
\begin{enumerate}
\item The variety $V_d$ is a Calabi--Yau hypersurface of dimension $d-2$ over $k(t)$. After the extension $t=s^d$, it is isomorphic to the base change of the Fermat hypersurface
\[
 F_d=Z(y_0^d+\cdots+y_{d-1}^d)\subset\PP^{d-1}_k.
\]
\item If $k$ is algebraically closed and $e\geq1$, then
\begin{equation}\label{equation:diagonal_spectrum}
 V_d\bigl(k\laurent{t^{1/e}}\bigr)\neq\emptyset
 \quad\Longleftrightarrow\quad \gcd(e,d)>1.
\end{equation}
The same assertion holds at $t=\infty$.
\end{enumerate}
\end{proposition}

\begin{proof} The first statement is straightforward: Adjunction gives $\omega_{V_d}\simeq\mathcal O_{V_d}$, the hypersurface
exact sequence gives $H^i(V_d,\mathcal O_{V_d})=0$ for $0<i<d-2$, and
the Lefschetz hyperplane theorem gives geometric simple connectedness
since $d-2\geq 2$. Moreover,  after $t=s^d$, the substitution $y_i=s^ix_i$ identifies $V_d$ with $F_d$. This proves~(i).

To prove~(ii),  we use the same valuation argument  as the proof of \cite[Proposition~4.4]{ELW}. Write $t=u^e$ and use the $u$-adic valuation. At a point of $V_d(k\laurent{u})$, the nonzero terms in the defining equation have valuations
\[
 ei+d\,\ord_u(x_i).
\]
If $\gcd(e,d)=1$, these valuations have distinct residues modulo $d$. There is therefore a unique term of least valuation, which is impossible in a sum equal to zero.

Conversely, suppose that $g=\gcd(e,d)>1$ and put $j=d/g$. Choose $\alpha\in k$ with $\alpha^d=-1$. The coordinates
\[
 x_0=1,\qquad x_j=\alpha u^{-e/g},\qquad x_i=0\quad(i\neq0,j)
\]
define a point, since $1+u^{ej}\alpha^d u^{-de/g}=0$.
At infinity, put $v=t^{-1}$, multiply the equation by $v^{d-1}$, and reverse the order of the coordinates. The result is again \eqref{equation:diagonal}, with parameter $v$.
\end{proof}

\begin{remark}
Suppose that $d=pq$ with $p<q$ primes, and let
$X\to\PP^1_{\CC}$ be a regular projective model of $V_d$.
In each of the fibres $X_0$ and $X_\infty$, the least component
multiplicity is $p$, and there is also a component of multiplicity $q$.
In particular, both fibres have gcd-multiplicity one.

By symmetry, it suffices to consider $X_0$. By
Proposition~\ref{proposition:diagonal}(ii), the positive integers $e$
for which $V_d(\CC\laurent{t^{1/e}})$ is nonempty are precisely the
multiples of $p$ or $q$. Lemma~\ref{lemma:local_index} therefore gives
$\operatorname{inf\text{-}mult}(X_0)=p$. To see that $q$ also occurs, write $X_0=\sum_j a_jF_j$.
By \cite[Lemma~3.1]{BCJW}, each $a_j$ is divisible by $p$
or $q$, and $q$ is a nonnegative integral linear combination of the
$a_j$. Since $p\nmid q$, some $a_j$ occurring in this combination
with positive coefficient is not divisible by $p$. Thus $q\mid a_j$,
and $a_j\leq q$ forces $a_j=q$.
\end{remark}

\begin{proof}[Proof of Theorem~\ref{theorem:B}]
Let $S=\sum_i D_i$, which is nonzero since $g(C)\leq 1$ and
$(C,\Delta)$ is of general type. We first make a base change of Lang's
family that leaves a single inf-multiple geometric fibre. We then pull
back to $C$ and place that fibre over the points of $S$.

Put $d=pq$, and let $V_d\subset\PP^{d-1}_{\QQ(t)}$ be Lang's diagonal
hypersurface, defined by $\sum_{i=0}^{d-1}t^ix_i^d=0$. Choose a regular
projective model $\mathcal V\to\PP^1_{\QQ}$ of $V_d$. Consider the
morphism
\[
  \phi\colon\PP^1_{\QQ}\longrightarrow\PP^1_{\QQ},
  \qquad z\longmapsto z^{q-p}(z-1)^p,
\]
which induces the field extension $\QQ(t)\hookrightarrow\QQ(z)$
given by $t\mapsto\phi(z)$. Let $\mathcal W\to\PP^1_{\QQ}$ be a regular
projective model of $V_d\otimes_{\QQ(t)}\QQ(z)$ obtained by resolving the pullback of $\mathcal{V}$ along $\phi$

The morphism $\mathcal{V}\to\mathbb{P}^1_{\mathbb{Q}}$ has two
inf-multiple fibres. However, as we show now, the morphism
$\mathcal{W}\to\mathbb{P}^1_{\mathbb{Q}}$ has a unique inf-multiple
fibre. Indeed, if $u$ is a local parameter and $t$ pulls back to a
function of order $a$, then this function is $u^a$ times a unit in
$\overline K[[u]]$. Every such unit is a $d$-th power, so the 
calculation in Proposition~\ref{proposition:diagonal}(ii) shows that
the pulled-back hypersurface has a point over
$\overline K\laurent{u^{1/e}}$ if and only if $\gcd(ea,d)>1$.
The orders of $\phi$ at $0,1,\infty$ are $q-p,p,-q$, respectively,
and its order is zero elsewhere. Thus at $z=0$ the degrees giving a
local point are precisely the multiples of $p$ or $q$, whereas at every
other point there is already a point over $\overline K\laurent{u}$.
By Lemma~\ref{lemma:local_index}, the base change therefore removes the
inf-multiplicity at infinity and leaves a single inf-multiple geometric
fibre, over $z=0$.

We now construct a suitable morphism from $C$ to $\mathbb{P}^1_K$.
Choose $n$ sufficiently large and coprime to $d$. Then, the linear system
$|nS|$ is basepoint-free and, since $K$ is infinite, contains a
$K$-divisor disjoint from $S$. Thus there is a rational function
$s\in K(C)$ whose zero divisor is exactly $nS$. Extend scalars to $K$
and pull back $\mathcal W_K\to\PP^1_K$ along $s\colon C\to\PP^1_K$.
Let $X$ be a projective desingularization of
$\mathcal W_K\times_{\PP^1_K}C$, so that the following diagram commutes
\[
\begin{tikzcd}[column sep=large]
X \arrow[r] \arrow[d,"f"'] &
\mathcal W_K \arrow[r] \arrow[d] &
\mathcal W \arrow[r] \arrow[d] &
\mathcal V \arrow[d] \\
C \arrow[r,"s"'] &
\PP^1_K \arrow[r] &
\PP^1_{\QQ} \arrow[r,"\phi"'] &
\PP^1_{\QQ}.
\end{tikzcd}
\]
Here $X$ is smooth and projective over $K$, and $f$ is flat and
surjective. Put $h=\phi(s)$ and denote the generic fibre of $f$ by
$V_h$. It is the base change of $V_d$ along $t\mapsto h$, so
\begin{equation}\label{equation:global_diagonal}
  V_h=Z\left(x_0^{pq}+h x_1^{pq}+\cdots+h^{pq-1}x_{pq-1}^{pq}\right)
  \subset \PP^{pq-1}_{K(C)}.
\end{equation}
By Proposition~\ref{proposition:diagonal}(i), $V_h$ becomes isomorphic
to the Fermat hypersurface $F_{pq}$ after adjoining a $pq$-th root of
$h$. Thus $f$ is isotrivial, proving~(i).

Let $c$ be a geometric point of $S$ and let $u$ be a local parameter
at $c$. Since $s$ has a zero of order $n$ at $c$ and $s-1$ is a unit
there, we have $\ord_c(h)=n(q-p)$. For every $e\geq 1$, the local
calculation above gives
\[
  V_h\bigl(\overline K\laurent{u^{1/e}}\bigr)\neq\emptyset
  \quad\Longleftrightarrow\quad
  \gcd\bigl(en(q-p),pq\bigr)>1
  \quad\Longleftrightarrow\quad
  p\mid e \ \text{or}\ q\mid e.
\]
The last equivalence uses  that $n$ is coprime to $pq$. In particular,  the minimal degree of a closed point on $V_{h, \overline K \laurent{u}}$  is $p$,
and there is also a point over the extension of degree $q$.
Lemma~\ref{lemma:local_index} therefore gives inf-multiplicity $p$
and gcd-multiplicity one. If $c\in C(\overline K)$ lies outside $S$,
then $s(c)\neq 0$. The local point at $s(c)$ pulls back to a point
over the completed function field at $c$, so the same lemma shows that
the fibre at $c$ has a component of multiplicity one. This proves~(ii).

By Proposition~\ref{proposition:diagonal}(i),   every smooth   fibre of $X\to C$ is a Calabi--Yau hypersurface. Thus, every such fibre is     special by \cite[Theorem~5.1]{Ca04},
hence weakly special. Since $g(C)\leq 1$, the curve $C$ is weakly
special. As $f$ has no divisible geometric fibres by~(ii),
\cite[Theorem~2.3]{BCJW} now implies that $X$ is weakly special.

Suppose that Orbifold Mordell holds for $(C,\Delta)$, and fix a finite
extension $L/K$. Since $\Delta_f=(1-1/p)S\geq\Delta$ by~(ii) and the
choice of $p$, the map $f\colon X\to(C,\Delta)$ is a morphism of
C-pairs. By properness and
\cite[Proposition~5.5]{BJCampanaStacks}\footnote{Campana credits Peyre
for the definition of orbifold rational points
\cite[Remark~9.22]{Ca04}. For the projection statement, see
\cite[Proposition~4.3]{Ca05} and Abramovich's formulation in terms of
soft integral points \cite[Proposition~2.1.8]{AbramovichBirGeom}.},
there is a finite set $T$ of finite places of $L$ such that every point
of $f(X(L))$ is a Campana $\mathcal O_{L,T}$-point of $(C,\Delta)$.
Orbifold Mordell therefore implies that $f(X(L))$ is finite.
\end{proof}
\begin{proof}[Proof of Theorem~\ref{theorem:A}]
Let $(C,\Delta)$ satisfy the hypotheses of the theorem. Choose primes
$p<q$ with $p$ at least as large as every multiplicity of $\Delta$,
and let $f\colon X\to C$ be the fibration of Theorem~\ref{theorem:B}.
By the last assertion of that theorem,  since $(C,\Delta)$ satisfies Orbifold Mordell (by assumption), the set $f(X(L))$ is finite for every
finite extension $L/K$.
Thus $X(L)$ lies in a finite union of fibres and is thus not Zariski dense
in $X$. Since $X$ is weakly special, this disproves the Weakly Special
Conjecture.
\end{proof}

\end{document}